# DILUTE ALGEBRAS AND SOLVABLE LATTICE MODELS

UWE GRIMM

*Instituut voor Theoretische Fysica, Universiteit van Amsterdam,*
*Valckenierstraat 65, 1018 XE Amsterdam, The Netherlands*

The definition of a dilute braid-monoid algebra is briefly reviewed. The construction of solvable vertex and interaction-round-a-face models built on representations of the dilute Temperley-Lieb and Birman-Wenzl-Murakami algebras is discussed.

## 1  Introduction

Solvable lattice models of two-dimensional statistical mechanics have been a major topic of research over the past decade. A most prominent rôle in this subject is played by the famous *Yang-Baxter equation* (YBE) representing a sufficient condition for solvability of a lattice model in the sense of commuting transfer matrices. A wealth of solutions to the YBE has been found; and an at least partial characterization of the solutions in terms of the classification of affine Lie algebras and their representations is available. The known solutions can be partitioned into two mutually dual classes, consisting of *vertex models* on the one hand and so-called *interaction-round-a-face* (IRF) models on the other. Fundamental solutions (in the sense that other solutions can be built from these by the fusion procedure) for vertex models related to non-exceptional affine Lie algebras have been constructed by Bazhanov[1] and by Jimbo,[2] corresponding IRF models for the non-twisted algebras by Jimbo, Miwa and Okado.[3] Later, Kuniba[4] also obtained solvable IRF models related to the twisted affine Lie algebras $A_n^{(2)}$.

Besides the interest in the physical properties of such systems, the wide attention attracted by the theory of solvable lattice models has been stimulated by the uncovery of numerous, partly surprising connections to other fields of mathematics and theoretical physics. Arguably, the most important implication in mathematics was the introduction of quantum groups, but also the relations to the theory of knot and link invariants has received considerable attention.[5–8] For the latter case, the connection can be established in an elegant way by considering so-called *braid-monoid algebras*, see e.g. the review by Wadati *et al.*[7] for details. The simplest example of a braid-monoid algebra is represented by the *Temperley-Lieb* (TL) algebra[9] corresponding to the $A_1^{(1)}$ models, i.e., the six-vertex and the ABF[10] models. However, these models can also be related to the *Birman-Wenzl-Murakami* (BWM) algebra,[11,12] which is the second example under consideration. In fact, the $B_n^{(1)}$, $C_n^{(1)}$, $D_n^{(1)}$ and $A_n^{(2)}$ models among the fundamental solutions of the YBE mentioned above are all related to the BWM algebra, while the $A_n^{(1)}$ models correspond to the Hecke algebra (respectively the TL algebra for $n = 1$). Of the non-exceptional series, solely the $D_n^{(2)}$ vertex models [1,2] do not seem to fit into this scheme.

The notion of a *dilute braid-monoid algebra* [13] naturally emerged after the discovery [14–17] of a second series of solvable IRF models related to $A_2^{(2)}$, the so-called dilute A–D–E models. These are quite different from Kuniba's models, and in particular include models which are solvable in the presence of a symmetry-breaking field. As these models can be described in terms of a *dilute Temperley-Lieb* (dTL) algebra, the investigation of the *dilute Birman-Wenzl-Murakami* (dBWM) algebra [18–21] and the related solvable lattice models (which include the $D_n^{(2)}$ vertex models [1,2]) was begun soon after. Recently, the more general case of a *two-colour BWM* algebra [13] has also been considered. [22]

In this note, we give a short account of the dilute case, commencing with a complete description of dilute braid-monoid algebras. In Sec. 3, the *Baxterization* [23] of the dTL and dBWM algebras are given. Subsequently, the corresponding vertex and IRF models are briefly discussed. Finally, we conclude with a short summary.

## 2  Dilute Braid-Monoid Algebras

A dilute braid-monoid algebra is an algebra with identity $I$ generated by the following three sets of generators

$$\mathcal{P}_j = \{()()_j, ()\,)_j, (\,()_j, (\,)_j\}, \qquad \mathcal{B}_j = \{(\times)_j, (\times)_j, (\underset{\frown}{\smile})_j\},$$
$$\mathcal{D}_j = \{(/)_j, (\backslash)_j, (\frown)_j, (\smile)_j\}, \tag{1}$$

where the subscript $j$ takes integer values $j = 1, 2, \ldots, N$. Here, we chose graphical symbols for the generators, anticipating the diagrammatical interpretation of the algebra. This is completely analogous to the case of a usual (non-dilute) braid-monoid algebra, which has a diagrammatic interpretation in terms of an array of $N+1$ *strings* or *strands*. The difference to the dilute case consists in the additional possibility that some of these strings are replaced by *vacancies* (empty strings).

Any generator with subscript $j$ acts non-trivially only at positions $j$ and $j+1$, which are the two strings respectively vacancies shown in our notation; and it is understood that the generators act on these two strings (vacancies) as indicated by their diagrams. [18,21] The generators in the set $\mathcal{P}_j$ are interpreted as *projectors* onto the four different choices of strings or vacancies at positions $j$ and $j+1$. For simplicity, we also introduce single-site projectors $s_j$ (for strings) and $v_j$ (for vacancies) by

$$\bigl(\rangle\langle\bigr)_j = s_j s_{j+1}, \quad \bigl(\rangle\ \bigr)_j = s_j v_{j+1}, \quad \bigl(\ \langle\bigr)_j = v_j s_{j+1}, \quad \bigl(\ \ \bigr)_j = v_j v_{j+1}. \quad (2)$$

These are required to fulfill the relations

$$s_j + v_j = I, \qquad s_j^2 = s_j, \qquad s_j s_k = s_k s_j \quad \text{for } j \neq k, \qquad (3)$$

which imply that $s_j$ and $v_j$, and thereby the generators in $\mathcal{P}_j$, are orthogonal projectors, i.e., $s_j v_j = v_j s_j = 0$.

This is consistent with the interpretation of the product of two generators as concatenation of the corresponding diagrams, where any diagram with a mismatch of strings or vacancies corresponds to the zero element of the algebra. In a similar way, all defining relations of the dilute braid-monoid algebra have a natural interpretation in terms of diagrams. To start with, there is a number of rather trivial commutation relations

$$O_j O_k = O_k O_j \qquad \text{for } |j-k| > 1 \qquad (4)$$

where $O_j \in \mathcal{P}_j \cup \mathcal{B}_j \cup \mathcal{D}_j$ and $O_k \in \mathcal{B}_k \cup \mathcal{D}_k$. Also, compatibility of the diagrams leads to

$$\begin{aligned}
\bigl(\rangle\langle\bigr)_j O_j \bigl(\rangle\langle\bigr)_j &= O_j \qquad \text{for all } O_j \in \mathcal{B}_j \\
\bigl(\ \langle\bigr)_j \bigl(/\bigr)_j \bigl(\rangle\ \bigr)_j &= \bigl(/\bigr)_j \qquad \bigl(\rangle\ \bigr)_j \bigl(\backslash\bigr)_j \bigl(\ \langle\bigr)_j = \bigl(\backslash\bigr)_j \\
\bigl(\rangle\langle\bigr)_j \bigl(\cup\bigr)_j \bigl(\ \ \bigr)_j &= \bigl(\cup\bigr)_j \qquad \bigl(\ \ \bigr)_j \bigl(\cap\bigr)_j \bigl(\rangle\langle\bigr)_j = \bigl(\cap\bigr)_j \qquad (5)
\end{aligned}$$

which one can think of as defining the 'external legs' of the corresponding generators.

The set $\mathcal{B}_j$, together with the $\bigl(\rangle\langle\bigr)_j$, generates a (non-dilute) braid-monoid subalgebra corresponding to all diagrams without vacancies. Therefore, the relations among these generators are just those of the usual braid-monoid algebra[7]

$$\begin{aligned}
(\times)_j (\times)_j &= (\times)_j (\times)_j = (\rangle\langle)_j \\
(\times)_j (\times)_{j+1} (\times)_j &= (\times)_{j+1} (\times)_j (\times)_{j+1} \\
(\asymp)_j^2 &= \sqrt{Q}\, (\asymp)_j \\
(\asymp)_j (\asymp)_{j\pm 1} (\asymp)_j &= (\asymp)_j (\rangle\langle)_{j\pm 1} \\
(\times)_j (\asymp)_j &= (\asymp)_j (\times)_j = \omega\, (\asymp)_j \\
(\times)_j (\times)_{j\pm 1} (\asymp)_j &= (\asymp)_{j\pm 1} (\times)_j (\times)_{j\pm 1} = (\asymp)_{j\pm 1} (\asymp)_j \qquad (6)
\end{aligned}$$

and can be interpreted graphically in the usual way as continuous deformations of the strings, where a closed loop yields a factor $\sqrt{Q}$ and removing a 'twist' gives rise to a factor $\omega$. The latter two objects are central elements in the algebra, hence we will always think of them as numbers. Note that the relations contain the defining relations of the braid group for the generators $(\diagdown\!\!\!\!\diagup)_j$ and $(\diagup\!\!\!\!\diagdown)_j$, and those of the TL algebra for the $(\asymp)_j$.

The remaining defining relations can be derived diagrammatically from those of Eq. 6 by simply replacing an arbitrary number of strings in the corresponding diagrams by vacancies.[18,21] A (certainly not unique) minimal set of relations obtained in this way reads

$$(/)_j (\backslash)_j = (\,()\,)_j \qquad (\backslash)_j (/)_j = ()\,)_j$$
$$(\cap)_j (\asymp)_j = \sqrt{Q}\,(\cap)_j \qquad (\asymp)_j (\cup)_j = \sqrt{Q}\,(\cup)_j$$
$$(\cap)_j (\cup)_j = \sqrt{Q}\,(\,)_j \qquad (\cup)_j (\cap)_j = (\asymp)_j$$
$$(/)_j (/)_{j+1} (\cup)_j = (\,()\,)_j (\cup)_{j+1} \qquad (\cap)_j (\backslash)_{j+1} (\backslash)_j = (\cap)_{j+1} (\,()\,)_j$$
$$(/)_{j+1} (/)_j = (\cap)_j (\cup)_{j+1} \qquad (\backslash)_j (\backslash)_{j+1} = (\cap)_{j+1} (\cup)_j$$
$$(/)_j (/)_{j+1} (\diagdown\!\!\!\!\diagup)_j = (\diagdown\!\!\!\!\diagup)_{j+1} (/)_j (/)_{j+1} \tag{7}$$

which completes the formal definition of a dilute braid-monoid algebra.

The two examples of interest are obtained as certain quotients of this algebra by demanding polynomial relations for the braid generators $(\diagdown\!\!\!\!\diagup)_j$, $(\diagup\!\!\!\!\diagdown)_j$ and the Temperley-Lieb generators $(\asymp)_j$, just as for the usual TL and BWM algebras.[7] For the dTL algebra, these relations have the form

$$0 = \left[(\diagdown\!\!\!\!\diagup)_j - q^{-1} (\,()\,)_j\right] \left[(\diagdown\!\!\!\!\diagup)_j + q^3 (\,()\,)_j\right], \qquad \omega = -q^3,$$
$$(\asymp)_j = q^{-1} \left[(\diagdown\!\!\!\!\diagup)_j - q^{-1} (\,()\,)_j\right], \qquad \sqrt{Q} = -(q^2 + q^{-2}); \tag{8}$$

and for the dBWM algebra they read as follows

$$0 = \left[(\diagdown\!\!\!\!\diagup)_j - q^{-1} (\,()\,)_j\right] \left[(\diagdown\!\!\!\!\diagup)_j + q (\,()\,)_j\right] \left[(\diagdown\!\!\!\!\diagup)_j - \omega (\,()\,)_j\right],$$
$$(\asymp)_j = \frac{\omega^{-1}}{q - q^{-1}} \left[(\diagdown\!\!\!\!\diagup)_j - q^{-1} (\,()\,)_j\right] \left[(\diagdown\!\!\!\!\diagup)_j + q (\,()\,)_j\right] = (\,()\,)_j + \frac{(\diagdown\!\!\!\!\diagup)_j - (\diagup\!\!\!\!\diagdown)_j}{q - q^{-1}},$$
$$\sqrt{Q} = 1 + \frac{\omega - \omega^{-1}}{q - q^{-1}}, \tag{9}$$

where $q = \exp(-i\lambda)$ is a complex number. Note that in both cases the relations for $\sqrt{Q}$ and $\omega$ follow from the other two equations.

## 3  Baxterizations

We now want to *Baxterize* [23,24] the dilute algebras. By this we mean the following. We want to find a general expression for a *local face operator* $X_j(u)$ ($u$ denoting the *spectral parameter*) which satisfy the YBE

$$X_{j+1}(u)\, X_j(u+v)\, X_{j+1}(v) \;=\; X_j(v)\, X_{j+1}(u+v)\, X_j(u) \qquad (10)$$

as well as the commutation relation

$$X_j(u)\, X_k(v) \;=\; X_k(v)\, X_j(u) \qquad \text{for } |j-k| > 1. \qquad (11)$$

as a consequence of the algebraic relations among the generators of the dilute algebra. In other words, the $X_j(u)$ satisfy the defining relations of a *Yang-Baxter algebra* (YBA).

For the dBWM algebra, the following expression

$$\begin{aligned}
X_j(u) \;=\;& \left(\rangle\langle\right)_j - \frac{\sin u}{2i \sin\lambda \, \sin\eta\lambda} \left[ e^{i(\eta\lambda - u)} \left(\diagup\!\!\!\!\diagdown\right)_j - e^{i(u - \eta\lambda)} \left(\diagdown\!\!\!\!\diagup\right)_j \right] \\
& + \frac{\sin u \, \sin(\eta\lambda - u)}{\sin\lambda \, \sin\eta\lambda} \left[ (\diagup)_j + (\diagdown)_j \right] + \frac{\sin(\eta\lambda - u)}{\sin\eta\lambda} \left[ (\rangle\,)_j + (\,\langle)_j \right] \\
& + \frac{\sin u}{\sin\eta\lambda} \left[ (\frown)_j + (\smile)_j \right] + \left[ 1 + \sigma \frac{\sin u \, \sin(\eta\lambda - u)}{\sin\lambda \, \sin\eta\lambda} \right] (\ \ )_j
\end{aligned} \qquad (12)$$

can be shown to satisfy Eqs. 10–11, where $\eta$ and $\sigma$ are defined by

$$q^{2\eta} \;=\; e^{-2i\eta\lambda} \;=\; \sigma\omega\,, \qquad \sigma^2 \;=\; 1\,. \qquad (13)$$

That is, from *any* representation of the dBWM algebra, characterized by certain values of $\eta$, $\sigma$, and possibly $\lambda$, one obtains a representation of the YBA via Eq. 12.

The case of the dTL algebra is very similar, but here $\eta$ and $\sigma$ take the *fixed* values $\eta = 3/2$ and $\sigma = -1$. Apart from that, the same expression as above Baxterizes the dTL algebra. Using Eq. 8, it can be rewritten in the slightly simpler form [13]

$$\begin{aligned}
X_j(u) \;=\;& \frac{\sin(\lambda - u)\,\sin(\eta\lambda - u)}{\sin\lambda \, \sin\eta\lambda} (\rangle\langle)_j - \frac{\sin u \, \sin(\eta\lambda - \lambda - u)}{\sin\lambda \, \sin\eta\lambda} (\asymp)_j \\
& + \frac{\sin u \, \sin(\eta\lambda - u)}{\sin\lambda \, \sin\eta\lambda} \left[ (\diagup)_j + (\diagdown)_j \right] + \frac{\sin(\eta\lambda - u)}{\sin\eta\lambda} \left[ (\rangle\,)_j + (\,\langle)_j \right] \\
& + \frac{\sin u}{\sin\eta\lambda} \left[ (\frown)_j + (\smile)_j \right] + \left[ 1 + \sigma \frac{\sin u \, \sin(\eta\lambda - u)}{\sin\lambda \, \sin\eta\lambda} \right] (\ \ )_j.
\end{aligned} \qquad (14)$$

by eliminating braids and inverse braids in favour of the TL generator $(\asymp)_j$.

## 4  Solvable Vertex and IRF Models

As the previous section shows, any representation of the dTL or the dBWM algebra give rise to a representation of the YBA, which in turn defines a solvable lattice model. Its Boltzmann weights are nothing but the matrix elements [7] of the local face operator $X_j(u)$; they are trigonometric functions of the spectral parameter $u$ because the relations defining the dilute braid-monoid algebras do not involve the spectral parameter. Depending on the type of representation under consideration, the resulting lattice model can be of vertex or IRF type, both with a finite or an infinite number of local states. Its Boltzmann weights are crossing symmetric [7] with crossing parameter $\eta\lambda$ and satisfy the local inversion relation

$$X_j(u)\,X_j(-u) \;=\; \varrho(u)\,\varrho(-u)\,I\,, \qquad \varrho(u) \;=\; \frac{\sin(\lambda-u)\,\sin(\eta\lambda-u)}{\sin\lambda\,\sin\eta\lambda}\,. \qquad (15)$$

Clearly, the next question is how to find suitable representations of the dilute algebras. In fact, this is rather simple, since one can make use of the many known representations of the TL and BWM algebra underlying the critical solvable vertex [1,2] and IRF [3] models of type $A_1^{(1)}$ (TL) and $B_n^{(1)}$, $C_n^{(1)}$ and $D_n^{(1)}$ (BWM), respectively. Any such representations can then be enlarged to a representation of the dilute algebra, where the original TL or BWM representation is contained as the restriction onto the completely occupied subspace (no vacancies).

In practice, this amount to enlarging the local space by one additional state [18,19] for representations of vertex type which live on an ($N$+1)-fold tensor product space. For IRF models, which are described in terms of adjacency diagrams restricting the allowed local configurations, the adjacency graph of the corresponding dilute model contains additional loops [21] that connect each local state to itself. In the case of IRF models of $B_n^{(1)}$ type, the original adjacency graphs already contain such loops which have to be distinguished from the additional dilute loops; basically, this amounts to introducing bond variables in the model.

Apart from this technical problem, the completion of the representation to that of the corresponding dilute algebra is straightforward. The only non-trivial matrix elements show up for the TL type generators $(\frown)_j$ and $(\smile)_j$. These can easily been found by realizing that the matrix elements of the TL generator $\left(\genfrac{}{}{0pt}{}{\smile}{\frown}\right)_j$ can be factored into products of two contributions as suggested by the diagram; and the corresponding matrix elements for the two dilute generators consists just of one of the two factors associated to the lower or upper part of the diagram.

Table 1: Representations of the TL and BWM algebras and their dilute counterparts.

| TL | $A_1^{(1)}$ $\rightsquigarrow$ $A_2^{(2)}$ | dTL |
|---|---|---|
| BWM | $\begin{cases} B_n^{(1)} & \rightsquigarrow & D_{n+1}^{(2)} \\ C_n^{(1)} & \rightsquigarrow & A_{2n}^{(2)} \\ D_n^{(1)} & \rightsquigarrow & B_n^{(1)} \end{cases}$ | dBWM |

As mentioned above, the solvable lattice models obtained in this way have trigonometric weights; they are critical if the underlying non-dilute model is critical. For many dilute IRF models, off-critical extensions involving elliptic functions have been found for both the dTL[15] and the dBWM[20,21] cases. In some cases, the elliptic nome breaks the $\mathbb{Z}_2$-symmetry of the adjacency diagram and hence acts as a symmetry breaking field rather than as a temperature-like variable which is the case in the underlying non-dilute models. With respect to their physical properties, only the dilute A–D–E models related to the dTL algebra have been studied in some detail,[15,16,25–29] and among those in particular the dilute $A_3$ model which can be regarded as a solvable version of a critical Ising model in a magnetic field.[15,16]

## 5 Concluding Remarks

We have presented an algebraic method to construct solvable lattice models. In some aspects, this method resembles the fusion procedure; however, the models we construct are completely different from those obtained by fusion. In general, the restricted IRF models (for which the value of the variable $\lambda$ is fixed by the adjacency diagram) obtained by *dilution* are new in the sense that they are not contained in the known series of non-dilute solvable lattice models and models obtained by fusion from those. However, the vertex models (and hence the corresponding $R$ matrices) are not. In this sense, the dilute models and the underlying representations of the dilute algebra are also labelled by affine Lie algebras as shown in Table 1. The corresponding similarity transformation for the $R$ matrices has been derived explicitly[18] for the $D_n^{(2)}$ models.

### Acknowledgments

The author is indebted to Bernard Nienhuis and Stichting FOM for giving him the opportunity to attend this conference. He also thanks the organizers M.-L. Ge and F. Y. Wu and the NANKAI Institute of Mathematics for their kind hospitality.